
\input epsf.tex
\input amssym.def
\input amssym
\magnification=1100
\baselineskip = 0.22truein
\lineskiplimit = 0.01truein
\lineskip = 0.01truein
\vsize = 8.5truein
\voffset = 0.2truein
\parskip = 0.10truein
\parindent = 0.3truein
\settabs 12 \columns
\hsize = 5.8truein
\hoffset = 0.2truein

\setbox\strutbox=\hbox{%
\vrule height .708\baselineskip
depth .292\baselineskip
width 0pt}
\font\caps=cmcsc10
\font\bigtenrm=cmr10 at 14pt

\def\sqr#1#2{{\vcenter{\vbox{\hrule height.#2pt
\hbox{\vrule width.#2pt height#1pt \kern#1pt
\vrule width.#2pt}
\hrule height.#2pt}}}}
\def\square{\mathchoice\sqr46\sqr46\sqr{3.1}6\sqr{2.3}4}

\centerline{\bigtenrm UNKNOTTING GENUS ONE KNOTS}
\tenrm
\vskip 14pt
\centerline{ALEXANDER COWARD AND MARC LACKENBY}
\vskip 18pt
\centerline{\caps 1. Introduction}
\vskip 6pt

There is no known algorithm for determining whether a knot has unknotting number
one, practical or otherwise. Indeed, there are many explicit knots ($11_{328}$ for example)
that are conjectured to have unknotting number two, but for which no proof of this
fact is currently available. For many years, the knot $8_{10}$ was in this class, but a
celebrated application of Heegaard Floer homology by Ozsv\'ath and Szab\'o [7] established
that its unknotting number is in fact two.

In this paper, we examine a related question: if a knot has unknotting number one,
are there only finitely many ways to unknot it, and if so, can one find them? This
is also a very difficult problem. However, we answer it completely here for knots
with genus one. 
We prove, in fact, that there is at most one way to unknot a genus one knot, with
the exception of the figure-eight knot, which admits precisely two unknotting
procedures.

We now make this statement more precise. To perform a crossing change to a knot $K$,
one proceeds as follows. A {\sl crossing circle} is a simple closed curve $C$ in the
complement of $K$ which bounds a disc that intersects $K$ transversely in two points
of opposite sign. If we perform $\pm 1$ surgery along $C$, then $K$ is transformed
into a new knot by {\sl changing a crossing}. A crossing change is
{\sl unknotting} if the resulting knot is the unknot. Two crossing changes are {\sl
equivalent} if the surgery coefficients are the same and there is an ambient isotopy, 
keeping $K$ fixed throughout, that takes one crossing circle to the other. 

\vskip 18pt
\centerline{
\epsfxsize=2in
\epsfbox{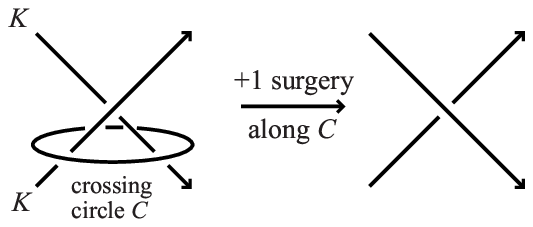}
}
\vskip 6pt
\centerline{Figure 1: a crossing change}

\noindent {\bf Theorem 1.1.} {\sl Suppose that $K$ is a knot with genus one and
unknotting number one. Then, if $K$ is not the figure-eight knot, there is
precisely one crossing change that unknots $K$, up to equivalence. If $K$ is the
figure-eight knot, then there are precisely two unknotting crossing changes, up to equivalence.}

It is possible to be rather more explicit about the crossing changes described above.
It is a theorem of Scharlemann and Thompson [9] that a knot has unknotting number one
and genus one if and only if it is a doubled knot. 
To construct the latter, one starts with the knot in the solid torus shown
in Figure 2, and then one embeds the solid torus into the 3-sphere.
The result is a {\sl doubled knot}, provided it is non-trivial.
If the solid torus is unknotted, then the resulting knot is known as a
{\sl twist knot}, again provided it is non-trivial.

\vskip 18pt
\centerline{
\epsfxsize=3.5in
\epsfbox{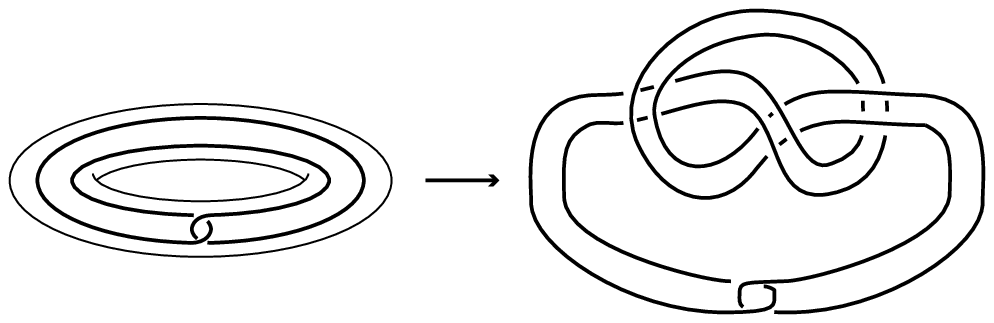}
}
\vskip 6pt
\centerline{Figure 2: a doubled knot}

\noindent {\bf Theorem 1.2.} {\sl If $K$ is a doubled knot, but not the figure-eight knot,
then the unique crossing circle that specifies an unknotting crossing change
is as shown in Figure 3. The two non-isotopic crossing circles that specify unknotting
crossing changes for the figure-eight knot are also shown in Figure 3.}

\vskip 18pt
\centerline{
\epsfxsize=3.1in
\epsfbox{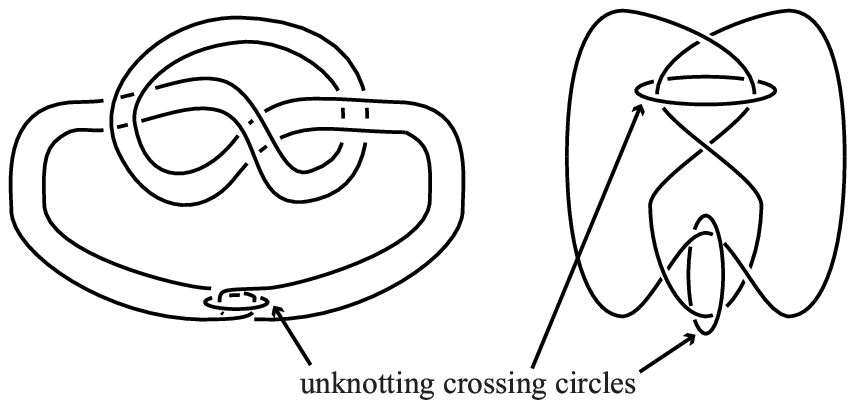}
}
\vskip 6pt
\centerline{Figure 3}

Theorem 1.2 is a rapid consequence of Theorem 1.1. This is because the crossing
circles in Figure 3 do indeed result in unknotting crossing changes. In the case
of the figure-eight knot, the two crossing circles can easily be verified to be
non-isotopic. There are several ways to prove this. One method is to show that they have
distinct geodesic representatives in the hyperbolic structure on the figure-eight knot
complement, and so are not even freely homotopic. Alternatively, observe that the
surgery coefficients are distinct for these two curves, and at most one choice of
surgery coefficient along a crossing circle can result in an unknotting crossing
change [8, Theorem 5.1]. Hence, these crossing circles cannot be ambient isotopic, via an isotopy
that fixes the figure-eight knot throughout.
Thus, by Theorem 1.1, these are the only unknotting crossing circles.

The case of the trefoil knot might, at first, be slightly worrying. It appears
to have three crossing circles that specify unknotting crossing changes, as shown
in Figure 4. However, these are all, in fact, equivalent. This can be seen directly.
It will also follow from the analysis in Section 6.

\vskip 18pt
\centerline{
\epsfxsize=1.8in
\epsfbox{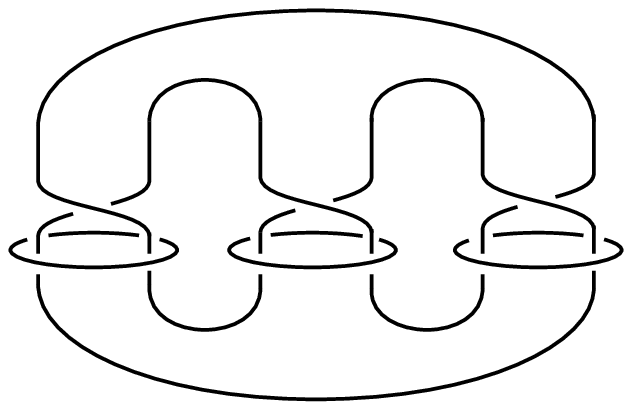}
}
\vskip 6pt
\centerline{Figure 4: unknotting crossing circles for the trefoil}

Theorems 1.1 and 1.2 are proved by extending some of the known techniques for detecting
unknotting number one knots. It is a theorem of Scharlemann and Thompson [9, Proposition 3.1]
that if a knot has unknotting number one, then it has a minimal genus Seifert surface that
is obtained by plumbing two surfaces, one of which is a Hopf band (see Figure
5). Moreover, the crossing change that unknots the knot has the effect of untwisting
this band. It is known that a knot $K$ which is not a satellite knot has, up to
ambient isotopy that fixes $K$, only finitely many minimal genus Seifert surfaces [12].
Thus, one is led to the following problem: can a minimal genus Seifert surface be obtained by plumbing
a Hopf band in infinitely many distinct ways, up to ambient isotopy that leaves the surface
invariant? Slightly surprisingly, the answer to this question is `yes'. However, it
is quite possible for different plumbings to result in the same
crossing change. The analysis for arbitrary knots quickly becomes difficult, but
when the knot has genus one, then progress can be made. In fact, for most doubled knots, 
it is fairly straightforward to show that there is a unique way
to unknot it. It turns out that the difficult case is when the knot is fibred. There
are precisely two fibred knots with genus one: the figure-eight and the trefoil. To
prove the theorem in these cases seems to require an analysis of the arc
complex of the once-punctured torus.

Given the finiteness result in Theorem 1.1, it is natural to enquire whether it
extends to higher genus knots. It seems reasonable to make the following conjecture.

\noindent {\bf Conjecture 1.3.} {\sl For any given knot $K$, there are only finitely
many crossing changes that unknot $K$, up to equivalence.}

\vfill\eject
Theorem 1.1 can be viewed as positive evidence for this. Further support comes from the
following result, of the second author (Theorem 1.1 of [3]). It deals with {\sl generalised crossing changes
of order $q$}, where $q \in {\Bbb N}$, which are defined to be the result of $\pm 1/q$ surgery along a crossing circle.

\noindent {\bf Theorem 1.4.} {\sl For any given knot $K$, there are only
finitely many generalised crossing changes of order $q$, where $q > 1$,
that unknot $K$, up to equivalence.}

However, there is still some reason to cautious about Conjecture 1.3. It can
be viewed as the knot-theoretic analogue to the following problem: does a given
closed orientable 3-manifold admit only finitely many descriptions as surgery
along a knot in 3-sphere? Examples due to Osoinach [6] show that the answer to
this question can be `no'.

\vskip 12pt
\centerline{
\epsfxsize=4in
\epsfbox{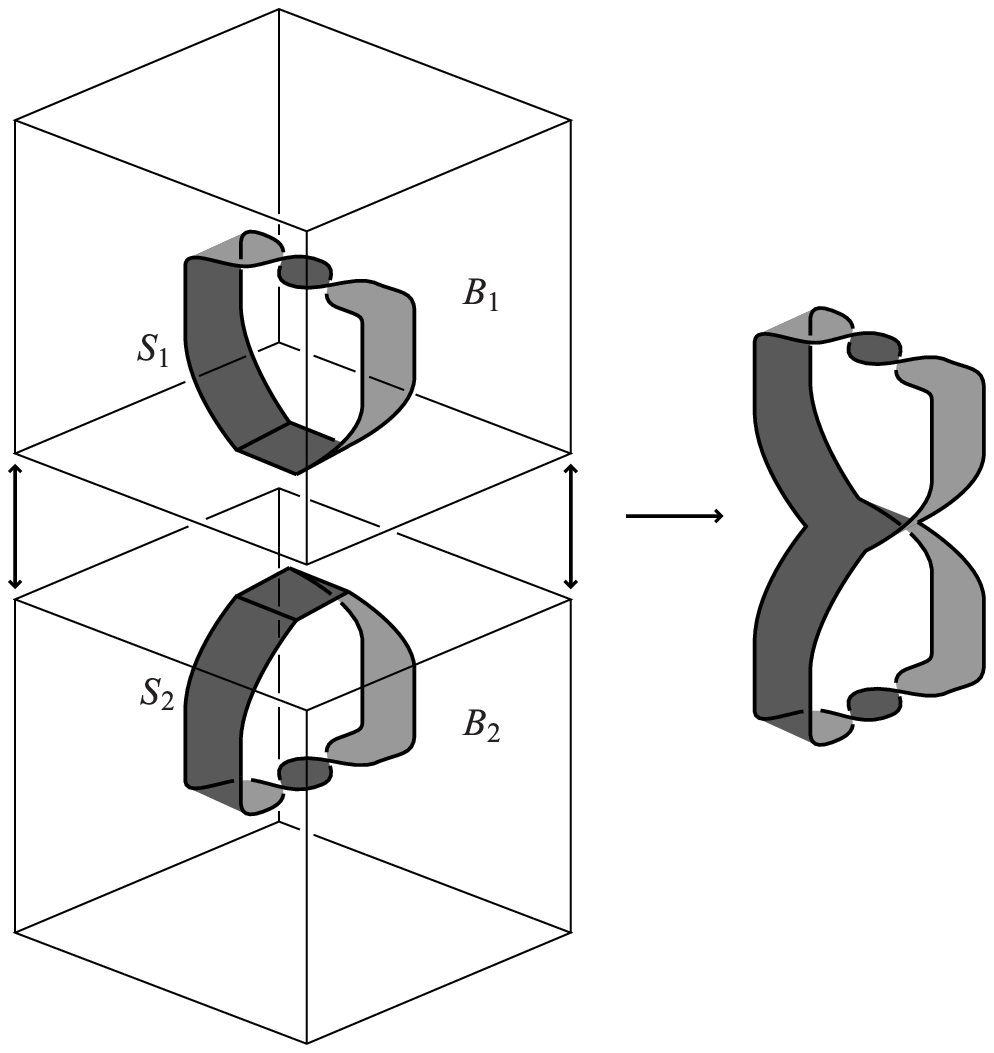}
}
\vskip 6pt
\centerline{Figure 5: plumbing}

\noindent {\sl Acknowledgements.} We would like to thank Peter Horn and Jessica Banks 
for pointing out errors in early versions of this paper. The second author was supported by
an EPSRC Advanced Research Fellowship.

\vfill\eject
\centerline{\caps 2. Plumbing and clean product discs}
\vskip 6pt

We now recall the operation of plumbing in detail. Suppose that $S_1$ and
$S_2$ are compact orientable surfaces embedded in 3-balls $B_1$ and $B_2$. Suppose that the
intersection of each $S_i$ with $\partial B_i$ is a square $I \times I$ such that
$(I \times I) \cap \partial S_1 = I \times \partial I$ and $(I \times I) \cap
\partial S_2 = \partial I \times I$. Then the surface in $S^3$ obtained by {\sl plumbing} $S_1$ and $S_2$
is constructed by gluing the boundaries of $B_1$ and $B_2$ so that
the two copies of $I \times I$ are identified in a way that preserves their product
structures. (See Figure 5.)

Suppose that a surface $S$ can be obtained by plumbing in two ways, by combining
$S_1 \subset B_1$ with $S_2 \subset B_2$, and by combining $S'_1 \subset B'_1$
with $S'_2 \subset B'_2$. We say that these are {\sl equivalent} if there exists
an ambient isotopy of the 3-sphere, leaving $S$ invariant throughout, that takes
$S_i$ to $S'_i$ ($i=1,2$) and $B_i$ to $B'_i$ ($i=1,2$).

Suppose that $S_1$ is a {\sl Hopf band}, which is an unknotted annulus embedded in $B_1$ with
a full twist. Then we will focus on some associated structure. The {\sl associated
crossing disc} $D$ is a disc embedded in the interior of $B_1$ which intersects $S_1$ in
a single essential arc in the interior of $D$. The boundary of this disc is the {\sl associated crossing circle}.
(See Figure 6.)

\vskip 6pt
\centerline{
\epsfxsize=2.1in
\epsfbox{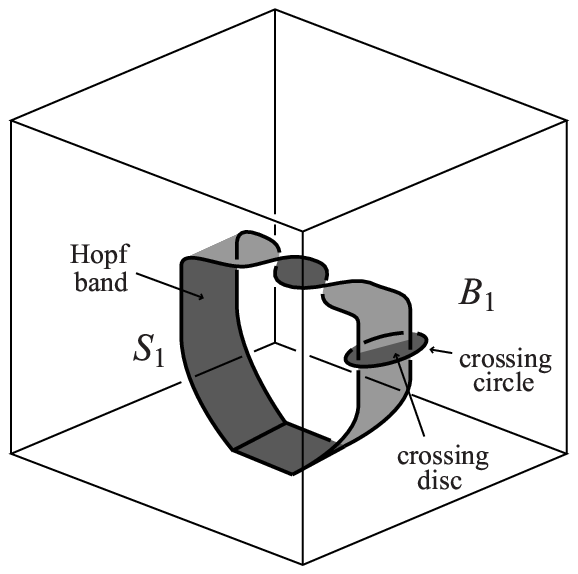}
}
\vskip 6pt
\centerline{Figure 6}

The following results are due to Scharlemann and Thompson [9, Proposition 3.1 and Corollary 3.2],
and are proved using sutured manifold theory [8].

\noindent {\bf Theorem 2.1.} {\sl Let $C$ be a crossing circle for a non-trivial knot $K$ such that
performing a crossing change along $C$ unknots $K$. Then $K$ has a minimal genus
Seifert surface which is obtained by plumbing surfaces $S_1$ and $S_2$,
where $S_1$ is a Hopf band. Moreover, there is an ambient isotopy, keeping $K$ fixed
throughout, that takes $C$ to the associated crossing circle for $S_1$.}

\noindent {\bf Corollary 2.2.} {\sl A knot has unknotting number one and genus one
if and only if it is a doubled knot.}

We now introduce some terminology.

Let $S$ be a Seifert surface for a knot $K$. Let $N(S)$ be a small product neighbourhood
$S \times I$. Let $S_-$ and $S_+$ denote the two components of $S \times \partial I$.

A {\sl product disc} is a disc properly embedded in $S^3 - {\rm int}(N(S))$ that
intersects $\partial S \times I$ in two vertical arcs. It therefore intersects $S_-$ in an arc
and $S_+$ in an arc. A {\sl direction} on a product disc $D$ is a choice of one
of the arcs $D \cap S_-$ or $D \cap S_+$. When such an arc is chosen, we say that
the product disc is {\sl directed towards} that arc, or just {\sl directed}.

A product disc $D$ is {\sl clean} if the projections of $S_- \cap D$ and $S_+ \cap
D$ to $S$ have disjoint interiors, up to ambient isotopy of $S^3 - {\rm int}(N(S))$ 
that fixes $\partial S \times I$ throughout.

Note that if $S$ is obtained by plumbing $S_1$ and $S_2$, where $S_1$ is a
Hopf band, then there is an associated directed clean product disc $D$ for $S$, defined as follows. The
intersection of $N(S)$ with $B_1$ is an unknotted solid torus. There is a unique
product disc embedded within $B_1$, as shown in Figure 7, which is clean and
which intersects $S_-$ and $S_+$ in essential arcs. We direct $D$ by
choosing the arc of $D \cap S_-$ and $D \cap S_+$ that avoids $S_2$.
Note that $D$ is {\sl essential}, in the sense that it forms a compression
disc for $\partial N(S)$ in $S^3 - {\rm int}(N(S))$. Also note that if $D$
is directed towards the arc $\alpha$, then the crossing circle associated
with the plumbing runs along $\alpha$ just above $S$, then
around $K$, then along $\alpha$ just below $S$ and then back around $K$.

\vskip 18pt
\centerline{
\epsfxsize=2.1in
\epsfbox{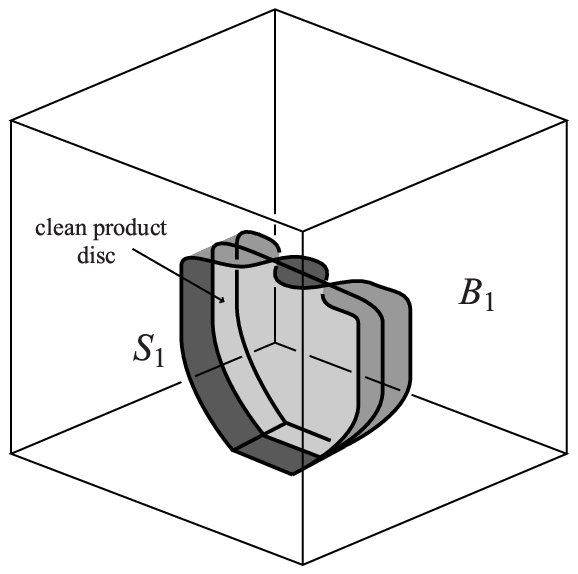}
}
\vskip 6pt
\centerline{Figure 7}

The projection of $D \cap S_-$ and $D \cap S_+$ to $S$ is two arcs $\alpha_-$ and
$\alpha_+$. We note that they have the following behaviour near $\partial S$. A
regular neighbourhood of $\alpha_+$ is a thickened arc, and we note that $\alpha_-$
intersects both sides of this thickened arc. Equivalently, $\alpha_+$ intersects
both sides of a regular neighbourhood of $\alpha_-$. When this is the case, 
we say that the arcs and the product disc
are {\sl alternating} (see Figure 8). Note that a clean alternating product disc is
automatically essential.

\vskip 18pt
\centerline{
\epsfxsize=1.8in
\epsfbox{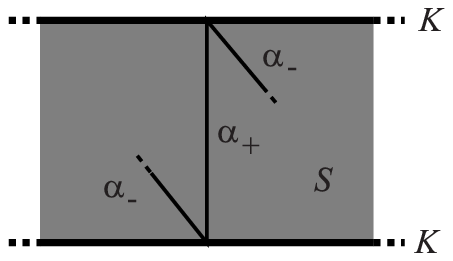}
}
\vskip 6pt
\centerline{Figure 8: alternating arcs}

The following is the main result of this section.

\noindent {\bf Theorem 2.3.} {\sl Let $S$ be a Seifert surface for a
knot $K$. Then there is a one-one correspondence between the following:
\item{(i)} decompositions of $S$, up to equivalence, as the plumbing of two surfaces, where the first surface 
is a Hopf band;
\item{(ii)} clean alternating directed product discs for $S$, up to ambient isotopy
that leaves $N(S)$ invariant and maintains the disc as a product disc throughout.

}

\noindent {\sl Proof.} We have already described how a plumbing of two surfaces as
in (i) determines a clean alternating directed product disc. Let us now determine how
a clean alternating directed product disc $D$ can be used to specify $S$ as the
plumbing of two surfaces, the first of which is a Hopf band.

Extend $D$ a little so that its boundary is $\alpha_- \cup \alpha_+$ in $S$.
Suppose that $D$ is directed towards $\alpha_+$.
Let $B$ be a small regular neighbourhood of $\alpha_+$ in
$S^3$. (See Figure 9.) We may ensure that $B$ intersects $S$ in a single disc $D_+$,
such that $D_+ \cap K$ is two properly embedded arcs in $B$. 
We may also ensure that $B \cap D$ is a single disc, and that
${\rm cl}(D - B)$ is a disc $D'$. Note that $D' \cap S$ is a
sub-arc of $\alpha_-$. Thicken $D'$ a little to form $D' \times I$.
We may ensure that $S \cap \partial B$ is an $I$-fibre in $D' \times I$,
because $\alpha_-$ and $\alpha_+$ are alternating. 
Then, $B \cup (D' \times I)$ is a 3-ball $B_1$. This will
be one of the 3-balls involved in the plumbing operation. The other
3-ball $B_2$ is the closure of the complement of $B_1$.

We now verify that this is indeed a plumbing construction.
We must check that $\partial B_1 \cap S$ is a disc $I \times I$
with the correct properties. By construction, $\partial B_1 \cap S$
is a disc. Identify this with $I \times I$, so that $\partial I \times I$ is
the intersection with $\partial B$. Then, near $\partial I \times I$, $S$ runs into
$B_1$, whereas near $I \times \partial I$, $S$ goes into $B_2$.
Thus, this is indeed a plumbing operation.

We must check that $B_1 \cap S$ is a Hopf band.
Note that $B_1 \cap S$ is $(I \times I) \cup D_+$, which is an 
annulus $A$. A regular neighbourhood of this annulus is a solid torus,
the boundary of which compresses in $B_1$ via a subdisc of $D$. Thus, $A$ is
unknotted. It has exactly one full twist because $D$ intersects
each component of $A \cap K$ in exactly one point (which is
one of the points of $\partial \alpha_+$). Thus, this is indeed
a plumbing of two surfaces as described in (i).

\vskip 6pt
\centerline{
\epsfxsize=2.8in
\epsfbox{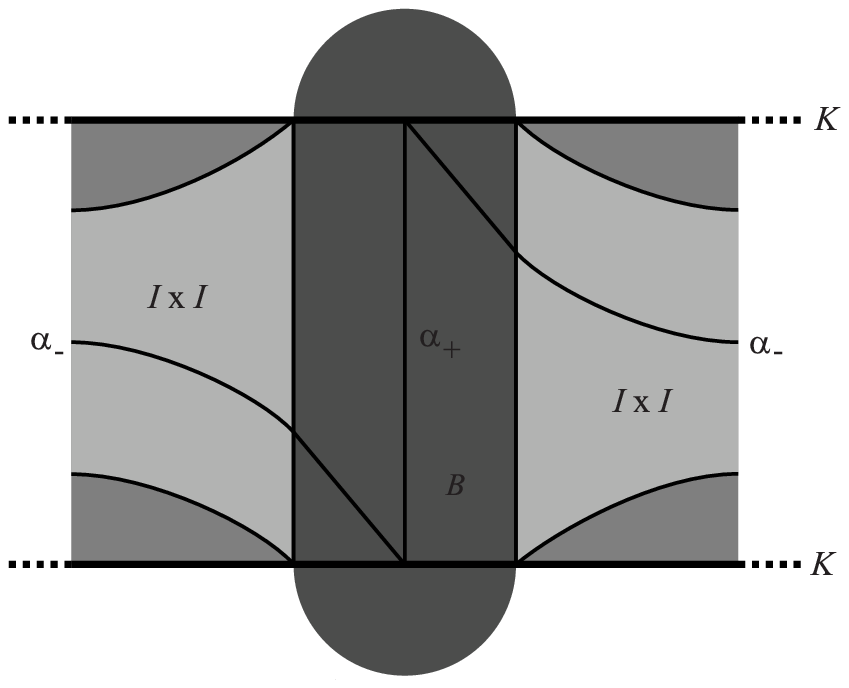}
}
\vskip 6pt
\centerline{Figure 9}

Finally, note that this does set up a one-one correspondence between
(i) and (ii). This is because an ambient isotopy of the clean product disc,
that leaves $S$ invariant, clearly results in an equivalent plumbing.
Conversely, equivalent plumbings result in isotopic clean product discs. 
$\square$

\noindent {\bf Remark 2.4.} This theorem implies that if a Seifert surface $S$
for a knot $K$ decomposes as a plumbing of two surfaces, where the first
surface is a Hopf band, then it does so in two ways which are inequivalent.
To see this, note that the plumbing determines a clean alternating directed
product disc $D$ by (i) $\rightarrow$ (ii) of the theorem. Now consider the
same disc $D$ but with the opposite direction. By (ii) $\rightarrow$ (i), we 
get a new way of describing $S$ as two plumbed surfaces. Associated with
these two ways of constructing $S$ via plumbing, there are two crossing circles.
These are actually ambient isotopic, via an isotopy that fixes $K$ throughout.
This is verified as follows. The first crossing circle runs along one
arc $\alpha_+$ of $S \cap D$ just above $S$, around $K$, then back along
$\alpha_+$ below $S$, then back around $K$. The second crossing circle
follows a similar route, but following the other arc $\alpha_-$ of $S \cap D$.
The arcs $\alpha_+$ and $\alpha_-$ are ambient isotopic, via an isotopy keeping $K$ fixed throughout,
since one may slide $\alpha_+$ across $D$ onto $\alpha_-$. This induces
an isotopy taking the first crossing circle to the second.

\vskip 18pt
\centerline{\caps 3. Fibred knots}
\vskip 6pt

When $K$ is a fibred knot, the following elementary observation allows us to translate
the existence of a clean product disc into information about the monodromy of the
fibration.

\noindent {\bf Proposition 3.1.} {\sl Let $K$ be a fibred knot with fibre $S$. Let $h \colon
S \rightarrow S$ be the monodromy, where $h | \partial S$ is the identity. Then,
there is a one-one correspondence between the following:
\item{(i)} clean essential product discs for $S$, up to ambient isotopy that leaves $N(S)$
invariant and maintains the disc as a product disc throughout;
\item{(ii)} properly embedded essential arcs $\alpha$ in $S$, up to isotopy of $\alpha$ in $S$, 
such that $h(\alpha)$ and $\alpha$ can be ambient isotoped, keeping their boundaries fixed,
so that their interiors are disjoint.

\noindent Moreover, the product disc is alternating if and only if $\alpha$ and $h(\alpha)$
are alternating.

}

\noindent {\sl Proof.} Since $K$ is fibred, the exterior of $S$ is a copy
of $S \times I$, where $S \times \partial I = S_- \cup S_+$. It is straightforward
that any essential product disc can be isotoped so that it respects the product structure
on $S \times I$. In other words, it is of the form $\alpha \times I$, for some
properly embedded essential arc $\alpha$ in $S$. Thus, $\alpha \times \{1\}$ is identified with $h(\alpha) \times \{0\}$.
If the disc is clean, then these arcs can be isotoped so that their interiors
are disjoint. Conversely, from any properly embedded essential arc $\alpha$ as in
(ii), one can construct a clean essential product disc from
$\alpha \times I$. $\square$


We note the following consequence of Proposition 3.1.

\noindent {\bf Corollary 3.2.} {\sl Let $K$ be a hyperbolic fibred knot with fibre $S$.
Suppose that $S$ is obtained by plumbing two surfaces, the first of which is a Hopf band.
Then it does so in infinitely many inequivalent ways.}

\noindent {\sl Proof.} Let $h$ be the monodromy of the fibration. Since $S$ is
obtained by plumbing two surfaces, the first of which is a
Hopf band, Theorem 2.3 states that it admits a clean alternating product
disc, which is automatically essential. Then, by Proposition 3.1, there is a properly embedded essential arc $\alpha$ in $S$
such that $h(\alpha)$ and $\alpha$ have disjoint interiors, up to isotopy in $S$, and are alternating. 
But now consider the arcs $h^n(\alpha)$ and
$h^{n+1}(\alpha)$, for each integer $n$. These also have disjoint interiors, up to isotopy in $S$, and are
alternating. Since $K$ is hyperbolic, the monodromy $h$ is
pseudo-Anosov and so the arcs $h^n(\alpha)$ are all distinct, up to isotopy. Thus, by Theorem 2.3
and Proposition 3.1, there are infinitely many distinct ways to decompose $S$ as the
plumbing of two surfaces, the first of which is a Hopf band. $\square$

As a consequence of the above proof, we make the following definitions.

Let $h \colon S \rightarrow S$ be a homeomorphism of a compact surface $S$.
We say that two arcs $\alpha_1$ and $\alpha_2$ on $S$ are {\sl $h$-equivalent} if
$h^n(\alpha_1)$ is ambient isotopic to $\alpha_2$ for some integer $n$.
This is clearly an equivalence relation.

Now let $K$ be a fibred knot with fibre $S$, and let $h$ be the monodromy.
We say that two arcs on $S$ are {\sl monodromy-equivalent}
if they are $h$-equivalent. 
We say that two clean essential directed product discs $D_1$ and $D_2$ are
{\sl monodromy-equivalent} if the directed product disc
$(h^n \times {\rm id})(D_1)$ is ambient isotopic to $D_2$, for some integer $n$, via
an isotopy that preserves direction. Here,
$h^n \times {\rm id}$ is a homeomorphism $S \times I \rightarrow S \times I$.
If these directed product discs are alternating, each one determines, by Theorem 2.3, a plumbing of $S$ into
two surfaces, the first of which is Hopf band. We also say that these
two plumbings are {\sl monodromy-equivalent}. Note that there is an ambient isotopy, leaving
$K$ fixed throughout, taking one plumbing decomposition to the other. This isotopy slides $S$
around the fibration. In particular, $S$ is not invariant. Nevertheless, we note that the associated crossing
circles are ambient isotopic, via an isotopy that leaves $K$ fixed throughout.

\vskip 18pt
\centerline{\caps 4. Seifert surfaces for doubled knots}
\vskip 6pt

According to Theorem 2.1, if a knot $K$ has unknotting number one, then some minimal genus
Seifert surface for $K$ is of a special form. If, in addition, $K$ has genus one, then
it is a doubled knot. It will therefore be useful to know the following.

\noindent {\bf Proposition 4.1.} {\sl Let $K$ be a doubled knot. Then $K$ has
a unique genus one Seifert surface, up to ambient isotopy that fixes $K$ throughout,
that is constructed by plumbing two surfaces, one of which is a Hopf band.}

Note that we are not claiming that the doubled knot $K$ has a unique
genus one Seifert surface. Indeed, this is not true in general, as Lyon showed in [5] that
a certain double of the $(3,4)$ torus knot has two distinct Seifert surfaces
with genus one.
However, by the above proposition, only one of these decomposes as the plumbing of
two surfaces, one of which is a Hopf band.

\noindent {\sl Proof.} Suppose that $K$ has a genus one Seifert surface $S$ that is
obtained by plumbing surfaces $S_1$ and $S_2$, where $S_1$ is a Hopf band.
Then $S_2$ is an annulus. Let $D$ be the clean product disc described in Section 2.
Then a regular neighbourhood of $N(S) \cup D$ is a solid torus $V$, with
boundary $T$ that lies in the complement of $K$. Inside $V$, the knot $K$ is as
shown in Figure 2. Note that $V - {\rm int}(N(K))$ is homeomorphic to the
exterior of the Whitehead link. This is hyperbolic, and so contains no essential
tori or annuli. We may assume that $V$
is embedded in the 3-sphere in a knotted fashion. Otherwise, $K$ is a 
twist knot and by [4], twist knots have a unique genus one Seifert surface,
up to ambient isotopy that fixes the knot throughout. Thus, $T$ is an
essential torus in the complement of $K$.

Suppose now that $K$ has another genus one Seifert surface $S'$ that is obtained by
plumbing two surfaces, one of which is a Hopf band. Then, we obtain another
solid torus $V'$ containing $K$ in the same way, with boundary $T'$. 
Since $T$ and $T'$ are incompressible in the complement of $K$ and since $V - {\rm int}(N(K))$ contains no 
essential annuli, there is an ambient isotopy of $T'$ in the complement of $K$ which renders it
disjoint from $T$. By switching the roles of $T$ and $T'$ if necessary, we may assume that $T'$ lies in $V$. Since
$V - {\rm int}(N(K))$ contains no essential tori, we deduce that $T'$ is
parallel to $T$. Hence, after a further small isotopy, we can
ensure that $V' = V$. Then $V$ contains both Seifert surfaces $S$ and $S'$.

The final step in the proof is to show that $V$ contains a unique genus one Seifert surface
for $K$, up to ambient isotopy that keeps $K$ fixed throughout. This is
proved in Sections 4 and 5 of [11], for example. Briefly, a proof
runs as follows. Let $D'$ be a compression disc for $T$ in $V$ that intersects
$S$ in a single embedded arc $\alpha$ in the interior of $D'$, running between the two points of $D' \cap K$.
We may pick $D'$ so that it is disjoint from the clean product disc $D$ for $S$.
By ambient isotoping $S'$ in $V$, keeping $K$ fixed, we may arrange that 
$D' \cap S'$ is a collection of properly embedded simple closed curves and
an arc $\alpha'$ in the interior of $D'$ running between the two points of $D' \cap K$.
We may also ensure that the intersection between $\alpha$ and $S'$
is $D' \cap K$. Now, by construction, $\alpha$ is a properly embedded essential
arc in $S$. Also, $\alpha'$ is properly embedded and essential in $S'$.
For if it were inessential, the resulting disc in $S'$ could be used to provide
an ambient isotopy of $K$ in $V$ that would render it disjoint from $D'$,
which is impossible. The exterior of $\alpha'$ in $S'$ is therefore an annulus,
since $S'$ has genus one. Any component of $(S \cap S') - K$ lies in this
annulus. We may arrange that $(S \cap S') - K$ is essential in $S$ and $S'$. Thus, if this intersection is non-empty,
it is a collection of parallel essential simple closed curves in $S$ and $S'$.
Since $S$ and $S'$ are homologous in $V - {\rm int}(N(K))$,
the intersection $S \cap S'$ is homologically trivial in $S$.
Hence, we can find adjacent curves of $S \cap S'$ in $S$ which
have opposite orientations, and which bound an annulus $A$. These are also 
parallel in $S'$, although not necessarily adjacent in $S'$,
and so we may ambient isotope the annulus between them in $S'$
onto $A$, and then eliminate these two intersection curves.
Thus, repeating this procedure if necessary, we may arrange that $S'$ is disjoint from $S$.

Now consider the clean product disc $D$ for $S$. We may arrange that $S' \cap D$
consists of a single arc that is essential in $S'$ and is not parallel in $S'$ to $\alpha'$. It therefore cuts
the annulus $S' - {\rm int}(N(\alpha'))$ into a disc. Perform an ambient isotopy of $S'$ that
takes $S' \cap D$ into $S$, and takes $\alpha'$ into $S$, but keeps 
the remainder of the interior of $S'$ disjoint from $S$. Hence,
$S' - S$ is the interior of a disc, and similarly $S - S'$ is the interior
of a disc. These discs are parallel in $V - K$, and
so there is an ambient isotopy of $V$, leaving $K$ fixed, taking $S'$ onto $S$. $\square$

\vskip 18pt
\centerline{\caps 5. Non-fibred genus one knots}
\vskip 6pt

\noindent {\bf Proposition 5.1.} {\sl Let $K$ be a non-fibred knot, with a genus one
Seifert surface $S$. Then, up to ambient isotopy of $S^3 - {\rm int}(N(S))$, it
admits at most one essential product disc.}

\noindent {\sl Proof.} Suppose that, on the contrary, there are two non-isotopic essential
product discs $D_1$ and $D_2$. We may ambient isotope these discs so that each component of 
$D_1 \cap D_2$ is an arc running from $S_-$ to $S_+$. Hence, a regular neighbourhood
of $(\partial S \times I) \cup D_1 \cup D_2$ in $S^3 - {\rm int}(N(S))$ 
is homeomorphic to $F \times I$, for some compact surface $F$, where
$F \times \partial I$ lies in $S_- \cup S_+$. If some component of
$\partial F \times \partial I$ bounds a disc in $S_-$ or $S_+$, enlarge $F \times I$
by attaching on a ball $D^2 \times I$. Repeating this as far as possible gives
an embedding of $F' \times I$ into $S^3 - {\rm int}(N(S))$, for some compact surface $F'$,
where $F' \times \partial I$ lies in $S_- \cup S_+$. Since $S$ is a once-punctured
torus, any two non-isotopic essential arcs in $S$ fill $S$, in the sense that the complement
of their union is a collection of discs. Hence, $F'$ is all of $S$ and so
$S^3 - {\rm int}(N(S))$ is homeomorphic to $S \times I$, the homeomorphism taking $S_- \cup S_+$
to $S \times \partial I$. This implies that $K$ is
fibred, which is contrary to hypothesis. $\square$

Combining this with Theorems 2.1 and 2.3, Remark 2.4 and Proposition 4.1, we deduce that there is at most one way
to unknot a non-fibred genus one knot. Thus, to complete the proof of
Theorem 1.1, we must now consider fibred genus one knots.

\vskip 18pt
\centerline{\caps 6. Fibred genus one knots}
\vskip 6pt

There are exactly two fibred genus one knots: the figure-eight knot and the trefoil \break [1, Proposition 5.14].
The monodromy of the figure-eight knot is pseudo-Anosov, whereas the monodromy of the
trefoil is periodic. Thus, the fibred case of Theorem 1.1 follows from Theorems 2.1 and 2.3,
Remark 2.4, Propositions 3.1 and 4.1 and the following.

\noindent {\bf Theorem 6.1.} {\sl Let $h \colon S \rightarrow S$ be an orientation-preserving homeomorphism
of the genus one surface with one boundary component, which is not isotopic to the identity. Then, up to $h$-equivalence,
there are at most two essential properly embedded arcs $\alpha$ in $S$
such that $\alpha$ and $h(\alpha)$ can be isotoped to be disjoint.
Moreover, if $h$ is periodic, then there is at most one such arc
up to $h$-equivalence.}

\vfill\eject
To prove this, we need to study the arc complex $C$
of $S$. Recall that this is defined as follows (see [2] for example).
It has a vertex for each ambient isotopy class of properly embedded essential unoriented arcs in $S$.
A collection of vertices span a simplex if and only if there are representatives
of the associated isotopy classes of arcs that are pairwise disjoint. The maximal
number of disjoint non-parallel essential arcs in $S$ is three, and hence the
complex is 2-dimensional.

There is a one-one correspondence between isotopy classes of properly embedded essential arcs on $S$
and ${\Bbb Q} \cup \{ \infty \}$. Thus, each properly embedded essential arc has an associated slope
$p/q$, and conversely each slope determines a unique isotopy class of properly embedded essential arc.
Two distinct slopes $p/q$ and $p'/q'$ (where these fractions are expressed in their
lowest terms) correspond to disjoint arcs if and only if $|pq' - qp'| = 1$.
Thus, the 1-skeleton of the arc complex is exactly the Farey graph.

Moreover, if one removes the vertices from the complex, the resulting space can be
identified with the hyperbolic plane ${\Bbb H}^2$ as follows. If we use the
upper-half space model for ${\Bbb H}^2$, its space at infinity $S^1_\infty$ is
${\Bbb R} \cup \{ \infty \}$. We place each vertex of the Farey graph
at the corresponding point of ${\Bbb Q} \cup \{ \infty \}$. We then
realise each open 1-simplex of $C$ as a geodesic. These geodesics divide ${\Bbb H}^2$
into ideal triangles. The closure of each such ideal triangle corresponds
to a 2-simplex of $C$. (See Figure 10.)

\vskip 18pt
\centerline{
\epsfxsize=2.5in
\epsfbox{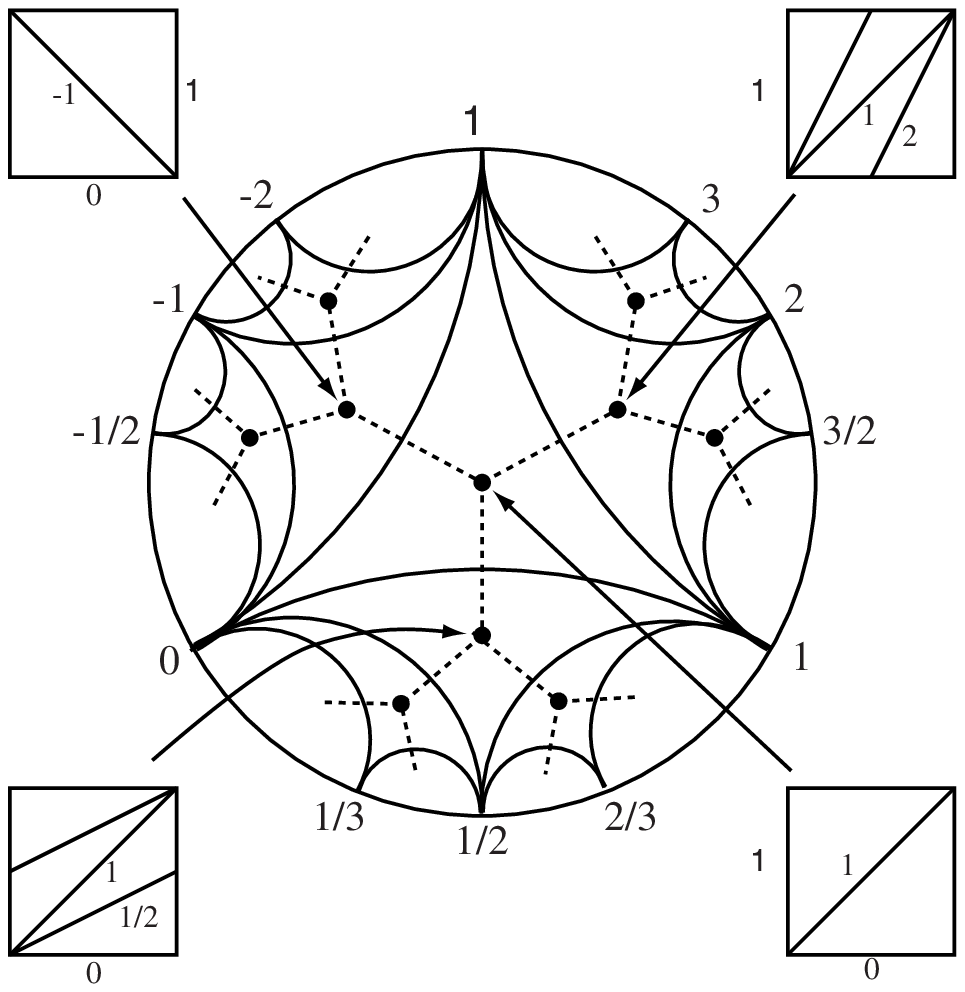}
}
\vskip 6pt
\centerline{Figure 10}

There is an associated dual complex $T$, which has a vertex at the centre of
each ideal triangle, and where two vertices are joined by an edge if and only if
their associated ideal triangles share an edge. This complex is a tree.

We will also need to collate some information about the link of a vertex $v$ of $C$.
By applying an automorphism to $C$, we may assume that $v = \infty$. The vertices
of $C$ adjacent to $\infty$ are precisely the integers. These are ordered around
$S^1_\infty - \{ \infty \}$, and successive integers are joined by an edge of $C$.
Thus, the link of $v$ in $C$ is a copy of the real line.

If $x_1, x_2, x_3, x_4$ are distinct points
in $S^1_\infty$ ordered that way around $S^1_\infty$, then we say that
$\{ x_1, x_3 \}$ is {\sl interleaved} with $\{ x_2, x_4 \}$.
Note that if the four points are vertices of the arc complex $C$, then
it cannot be the case both that $x_1$ and $x_3$ are adjacent in $C$,
and that $x_2$ and $x_4$ are adjacent in $C$. This is because the
corresponding edges of $C$ would then intersect at a point in ${\Bbb H}^2$,
which is impossible.

Any orientation-preserving homeomorphism $h$ of $S$ induces an automorphism of $H_1(S)$, and
hence is an element of ${\rm SL}(2,{\Bbb Z})$. The image of this in
${\rm PSL}(2, {\Bbb Z})$ corresponds to a M\"obius map, which
induces an orientation-preserving isometry of ${\Bbb H}^2$. There are three types of
orientation-preserving hyperbolic isometry: elliptic, parabolic and
loxodromic. These correspond to the three types of surface automorphism:
periodic, reducible (with infinite order) and pseudo-Anosov.

The isometry of ${\Bbb H}^2$ leaves $T$ invariant and induces a simplicial
automorphism of $C$. The automorphism of $T$ either fixes a point or is fixed-point free.
In the latter case, Proposition 24 in [10] implies that the isometry of $T$ has a unique
invariant line, known as the {\sl axis} of the automorphism. 
The automorphism acts by translation along this line, and therefore has
infinite order. The case where there is a fixed point in $T$ occurs exactly
when the hyperbolic isometry is elliptic and hence is periodic.

\noindent {\sl Case 1.} $h$ is periodic.

The induced automorphism of $T$ therefore fixes a 
point in $T$, which is either a vertex or a midpoint of an edge.
Let us consider first where the fixed point is a vertex of $T$.
Then, dually, there is an invariant 2-simplex in $C$. Thus, $h$ induces
a rotation of ${\Bbb H}^2$ about the centre of the corresponding ideal
triangle $\Delta$. We claim that the only properly embedded essential
arcs $\alpha$ on $S$ such that $h(\alpha) \cap \alpha = \emptyset$, up to isotopy, correspond
to ideal vertices of $\Delta$. To see this, note that, for any point $x$ on
$S^1_\infty$ that is not an ideal vertex of $\Delta$, $\{ h(x), x \}$ is
interleaved with two ideal vertices of $\Delta$. Since the latter form the
endpoints of an edge of $C$, we deduce that $h(x)$ and $x$ do not.
This proves the claim. Note now that the three vertices of $\Delta$ are
all $h$-equivalent. Hence, Theorem 6.1 is proved in this case.
Suppose now that $h$ fixes a midpoint of an edge $e$ of $T$ that
is dual to an edge $\overline e$ of $C$.
It therefore acts as a rotation about this point.
Hence,  the only properly embedded essential
arcs $\alpha$ on $S$ such that $h(\alpha) \cap \alpha = \emptyset$, up to isotopy, correspond
to endpoints of $\overline e$. This is because, for any point $x \in S^1_\infty - \partial \overline e$,
$\{ x, h(x) \}$ is interleaved with $\partial \overline e$. Thus, $x$ and $h(x)$
cannot be joined by an edge of $C$. Since the two endpoints of $\overline e$ are
$h$-equivalent, Theorem 6.1 is proved in the periodic case, and so there
is a unique unknotting crossing change for the trefoil knot.

\noindent {\sl Case 2.} $h$ is reducible and has infinite order.

Then $h$ leaves an essential arc invariant up to isotopy. Thus, the element of ${\rm SL}(2, {\Bbb Z})$ is
conjugate to $$\left ( \matrix { \pm 1 & n \cr 0 & \pm 1 } \right ),$$
where $n \in {\Bbb Z} - \{ 0 \}$.
If we can establish Theorem 6.1 for an automorphism $h \colon S \rightarrow S$,
then it also holds for any conjugate automorphism. Thus, we may assume
that $h$ is given by the above matrix.
It is then easy to see that, up to $h$-equivalence, the only properly embedded essential arcs
$\alpha$ in $S$ such that $h(\alpha)$ and $\alpha$ can be made disjoint
are $\infty$ and $0$, and that the latter only arises
if $n = \pm 1$. In particular, Theorem 6.1 holds in this case.

\noindent {\sl Case 3.} $h$ is pseudo-Anosov.

Then the induced action on ${\Bbb H}^2$ is loxodromic. The fixed points on $S^1_\infty$
cannot be vertices of $C$, because $h$ is not reducible. These are the endpoints
$A_-$ and $A_+$ of the axis $A$ on $S^1_\infty$. We may suppose that $A_-$ (respectively,
$A_+$) is the repelling (respectively, attracting) fixed point.

We say that two distinct vertices of $C$ are
{\sl on the same side of $A$} if they are not interleaved with the endpoints
of $A$. If $x_1$ and $x_2$
are distinct points on the same side of $A$, then we say that $x_1 < x_2$ if
$\{ x_1, A_+ \}$ and $\{ x_2, A_- \}$ are interleaved. This is a total
order on points on one side of $A$, which $h$ preserves.

We say that a vertex $v$ of $C$ is {\sl visible} from some point
$x \in A$, and that $x$ is {\sl visible} from $v$, if $x$ lies on an edge of $T$ that is dual to
an edge of $C$ that is incident to $v$. We say that $v$ is
{\sl visible from $A$} if it is visible from some point
of $A$. An equivalent condition is that some edge
emanating from $v$ is dual to an edge of $A$. Another equivalent condition is
that $v$ is adjacent in $C$ to a vertex that lies on the opposite side
of $A$.

\noindent {\sl Claim.} Suppose that, for some vertex $v$ of $C$,
$h(v)$ and $v$ are adjacent in $C$. Then, $v$ is visible from $A$,
and moreover, any vertex of $C$ that is visible from $A$ and that 
lies on the same side of $A$ as $v$ is $h$-equivalent to $v$.

The claim rapidly implies the theorem. This is because $A$ has only two sides,
and so there are at most two $h$-equivalence classes of vertices $v$
such that $h(v)$ and $v$ are adjacent in $C$. 

To prove the first part of the claim, suppose that $v$ is not visible
from $A$. The endpoints of the edges of $C$ emanating from $v$ are
arranged in order around $S^1_\infty$, and successive vertices in this
ordering are adjacent in $C$. Since $v$ is not visible, there are two
successive vertices $v_-$ and $v_+$ such that a cyclic ordering of
vertices around $S^1_\infty$ is $v, v_-, A_-, A_+, v_+$. In other words,
$v_- < v < v_+$. Now, $h(v_+) > v_+$. Since $\{v_-, v_+\}$ and $\{ h(v_-), h(v_+) \}$
both form the endpoints of edges, they are not interleaved. Hence, 
$h(v_-) \geq v_+$. Since $v>v_-$, $h(v) > h(v_-)$ and so 
$h(v) > v_+$. Hence, $\{ h(v), v \}$ and $\{ v_-, v_+ \}$ are interleaved,
and therefore $h(v)$ and $v$ are not adjacent in $C$. This proves the
first part of the claim.

To prove the second part, we will construct a fundamental domain
for the action of $h$ on $A$, and will show that the only vertices
on the $v$ side of $A$ that are visible from this fundamental domain
are $v$ and $h(v)$. Hence, the only vertices that are visible from $A$
and that lie on the $v$ side of $A$ are $h$-equivalent to $v$.

By the first part of the claim, there is at least one vertex adjacent
to $v$ that does not lie on the same side of $A$ as $v$. There are only
finitely many such vertices, and so there is one $v'$ that is maximal.
Let $v_+$ be the vertex incident to both $v$ and $v'$, and that
lies on the $v$ side of $A$, and which satisfies $v < v_+$. Now,
$h(v)$ cannot be greater than $v_+$, for otherwise $\{ h(v),v \}$ and
$\{ v_+, v' \}$ would be interleaved. Also, $h(v)$ cannot be less than
$v_+$, since then $\{ h(v), h(v') \}$ and $\{ v, v_+ \}$ would be
interleaved. Hence, $h(v) = v_+$. 

\vskip 6pt
\centerline{
\epsfxsize=2.2in
\epsfbox{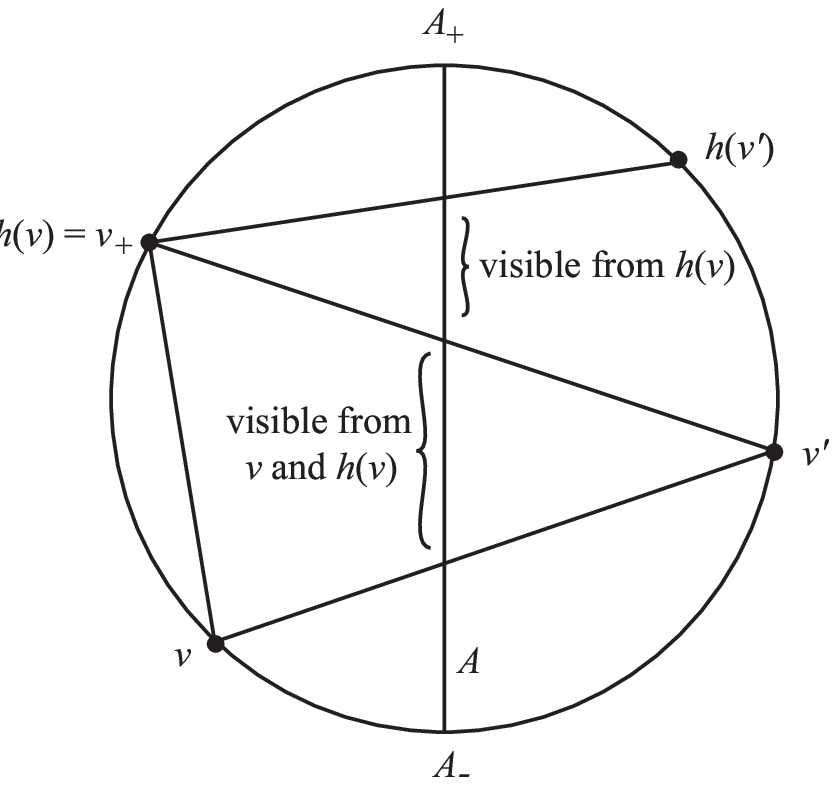}
}
\vskip 6pt
\centerline{Figure 11}

The required fundamental domain in $A$ is the interval that lies between
the edge $\{v,v'\}$ and the edge $\{ h(v), h(v') \}$.
This interval is divided into two sub-intervals by the edge
$\{ h(v), v' \}$. For each point in the first sub-interval, the only visible
vertices on the $v$ side of $A$ are $v$ and $h(v)$. For each point in the second sub-interval,
the only visible vertex on the $v$ side of $A$ is $h(v)$. This proves the claim, and hence the theorem.
$\square$

Thus, we have proved that there are at most two ways to unknot the
figure-eight knot. This is the last step in the proof of Theorem 1.1.
$\square$

\vskip 18pt
\centerline {\caps References}
\vskip 6pt

\item{1.} {\caps G. Burde, H. Zieschang}, 
{\sl Knots}, de Gruyter Studies in Mathematics, 5. Walter de Gruyter \& Co., Berlin, 1985.

\item{2.} {\caps A. Hatcher, W. Thurston,} {\sl Incompressible surfaces in 2-bridge knot
complements}, Invent. Math. 79 (1985) 225--246.

\item{3.} {\caps M. Lackenby,} {\sl Exceptional surgery curves in triangulated 
3-manifolds,}  Pacific J. Math.  210  (2003) 101--163.

\item{4.} {\caps H. Lyon}, {\sl Simple knots with unique spanning surfaces},
Topology 13 (1974) 275--279.

\item{5.} {\caps H. Lyon}, {\sl Simple knots without unique minimal surfaces},
Proc. Amer. Math. Soc. 43 (1974) 449--454.

\item{6.} {\caps J. Osoinach,} {\sl Manifolds obtained by surgery on an infinite
 number of knots in $S^3$,} Topology 45 (2006) 725--733.

\item{7.} {\caps P. Ozsv\'ath, Z. Szab\'o,}
{\sl Knots with unknotting number one and Heegaard Floer homology.}
Topology 44 (2005) 705--745.

\item{8.} {\caps M. Scharlemann,} {\sl Sutured manifolds and generalized Thurston norms,}
J. Differential Geom. 29 (1989) 557--614.

\item{9.} {\caps M. Scharlemann, A. Thompson,} {\sl Link genus and the Conway moves,} 
Comment. Math. Helv.  64 (1989) 527--535. 

\item{10.} {\caps J-P. Serre}, {\sl Arbres, Amalgames, ${SL}_2$}, Ast\'erisque 46,
Soc. Math. France (1977).

\item{11.} {\caps W. Whitten,} {\sl Isotopy types of knot spanning surfaces},
Topology 12 (1973) 373--380.

\item{12.} {\caps R. Wilson,} {\sl Knots with infinitely many incompressible Seifert surfaces},
J. Knot Theory Ramifications 17 (2008) 537--551.

\vskip 12pt
\+ Mathematics Department, University of California at Davis, \cr
\+ One Shields Avenue, Davis, California 95616, USA \cr

\vskip 12pt
\+ Mathematical Institute, University of Oxford, \cr
\+ 24-29 St Giles', Oxford OX1 3LB, United Kingdom. \cr

\end